\documentclass[12pt]{amsart}

\usepackage{amsmath}
\usepackage{latexsym}
\usepackage{amssymb}
\usepackage{amsfonts}

\include{PDF}
\usepackage{graphics}
\renewcommand{\proof}{\par\noindent{\it Proof.\ \ }}
\def\qed{\ifmmode\square\else\nolinebreak\hfill
$\square$\fi\par\vskip12pt}

\def\l{\langle} \def\r{\rangle}

\def\FF{\mathbb F} \def\ZZ{\mathbb Z}

\def\Cay{{\sf Cay}}

\def\calO{{\mathcal O}}

\def\D{{\rm D}} \def\Q{{\rm Q}}
\def\S{{\rm S}} \def\G{{\rm G}}
 \def\M{{\rm M}}
\def\soc{{\sf soc}}

\def\Z{{\bf Z}}

\def\Aut{{\sf Aut}}  \def\Out{{\sf Out}}
\def\K{{\sf K}}

\def\calO{{\mathcal O}}

\def\Ome{{\it\Omega}}
\def\Ga{{\it\Gamma}}

\def\SigmaL{{\rm \Sigma L}}

   \def\s{\sigma}
  \def\o{\omega}

\def\POmega{{\rm P\Omega}}
\def\Sp{{\rm Sp}}
\def\GammaL{{\rm \Gamma L}}

\def\PGammaL{{\rm P\Gamma L}} 

\def\A{{\rm A}}

\def\PSL{{\rm PSL}}
\def\GL{{\rm GL}} \def\SL{{\rm SL}}

 \def\PSU{{\rm PSU}}  
\def\Sz{{\rm Sz}}
\def\Ree{{\rm Ree}}

\def\Ra{{\sf rad}}

\def\McL{{\rm McL}} 
 
  \def\D{{\rm D}} \def\G{{\rm G}}
 
\def\Co{{\rm Co}} 
\def\HS{{\rm HS}}
\renewcommand{\leq}{\leqslant}
\renewcommand{\geq}{\geqslant}

\newtheorem{theorem}{Theorem}[section]%
\newtheorem{lemma}[theorem]{Lemma}%
\newtheorem{corollary}[theorem]{Corollary}%
\newtheorem{proposition}[theorem]{Proposition}%
\newtheorem{example}[theorem]{Example}%

\begin{document}

\title[Locally $2$-transitive graphs]
{On the Stabilisers of  Locally $2$-Transitive Graphs}
\thanks{{\it 2010 Mathematics subject classification}: 05E18, 20B25}
\thanks{This work forms a part of an ARC grant project.}

%
%

\author[Song]{Shu Jiao Song}
\address{School of Mathematics and Statistics\\
The University of Western Australia\\
Crawley, WA 6009, Australia}

\email{shu-jiao.song@uwa.edu.au}

\date\today

\maketitle

\begin{abstract}
For a connected locally $(G,s)$-arc-transitive graph $\Ga$ with $s\geqslant 2$ and an edge $\{v,w\}$,
determining the amalgam $(G_v,G_w,G_{vw})$ is a fundamental problem in the area of 
symmetrical graph theory, but it is very difficult.
In this paper, we give a classification of $(G_v,G_w,G_{vw})$ in the case where 
the vertex stabilisers $G_v$ and $G_w$ are faithful on their neighbourhoods,
which shows that except for the case $G_v\cong G_w$, there are exactly 16 such triples.
\end{abstract}

\section{Introduction}\label{Intro}

Let $\Ga=(V,E)$ be a connected undirected simple graph.
An {\it $s$-arc} of $\Ga$ is an $(s+1)$-tuple $(v_0,v_1,\dots,v_s)$ of
vertices  such that $\{v_{i-1},v_i\}\in E$ for $1\le i\le s$ and
$v_{i-1}\not=v_{i+1}$ for $1\le i\le s-1$.
For a group $G\le\Aut\Ga$, $\Ga$ is called {\it locally $(G,s)$-arc-transitive} if, for any vertex $v\in V$, the vertex stabiliser
$G_v$ acts transitively on the set of $t$-arcs starting at $v$, for any $t$ with $1\leq t\leq s$.
A locally $(G,s)$-arc transitive graph is called {\it $(G,s)$-arc transitive} if $G$ is also transitive on the vertex set.
Let $\Ga(v)$ be the set of vertices adjacent to $v$, and 
$G_v^{\Ga(v)}$ the permutation group on $\Ga(v)$ induced by $G_v$.

For a connected locally $(G,s)$-arc-transitive graph $\Ga$ and an edge $\{v,w\}$, 
the triple $(G_v,G_w,G_{vw})$ is called the {\it amalgam} of $G$, and also of $\Ga$.
A fundamental problem for studying locally $s$-arc-transitive graphs is to determine their amalgams.
van Bon \cite{vanBon-2} proved that if $G_{vw}^{[1]}=1$ then $\Ga$ is not locally $(G,s)$-arc-transitive.
Potocnik \cite{Potocnik} determined the amalgams for valency $\{3,4\}$ in the case where $G_{vw}^{[1]}=1$.
This paper is one of a series of papers which aim to classify the amalgams $(G_v,G_w,G_{vw})$ with
trivial edge kernel.

Assume that $\Ga$ is locally $(G,s)$-arc-transitive with $s\ge2$.
For a vertex $v$ of $\Ga$, denote by $G_v^{[1]}$ the kernel of $G_v$ acting on $\Ga(v)$.
Then $G_v^{\Ga(v)}\cong G_v/G_v^{[1]}$, and $G_v^{\Ga(v)}$ is a 2-transitive permutation group.
As usual, denote by $X^{(\infty)}$ the smallest normal subgroup of $X$ such that $X/X^{(\infty)}$ is soluble.
Let $\Ra(H)$ be the soluble radical of a group $H$, that is, the largest soluble normal subgroup of $H$.

Here we generalise a classical result for primitive permutation groups, see \cite[Section 18]{Wielandt}, 
which shows that some information of $(G_v,G_w,G_{vw})$ can be obtained from the permutation groups 
$G_v^{\Ga(v)}$ and $G_w^{\Ga(w)}$.

\begin{theorem}\label{key-lem}
Let $\Ga$ be a connected $G$-edge-transitive graph, where $G\leqslant\Aut\Ga$.
Then, for an edge $\{v,w\}$ of $\Ga$, the following statements hold.
\begin{itemize}
\item[(1)] Each composition factor of $G_v$ is a composition factor of $G_v^{\Ga(v)}$, $G_{vw}^{\Ga(w)}$ or 
$G_{vw}^{\Ga(v)}$.

\item[(2)] Each composition factor of $G_{vw}$ is a composition factor of $G_{vw}^{\Ga(w)}$ or $G_{vw}^{\Ga(v)}$.
\end{itemize}
\end{theorem}

\begin{theorem}\label{faith-stab}
Let $\Ga$ be a connected locally $(G,2)$-arc-transitive graph, and let $\{v,w\}$ be an edge of $\Ga$.
Assume that $G_v^{[1]}=G_w^{[1]}=1$.
Then either $G_v\cong G_w$ and $\Ga$ is regular, or $(G_v,G_w,G_{vw})$ is one of the sixteen triples listed in Table~$\ref{Intro}$.
\end{theorem}

\[\begin{array}{|c ccccccc |}\hline
&G_v&G_w&G_{vw}&|\Ga(v)|&|\Ga(w)|&G_{w_{-1}vw}&G_{vwv_1} \\ \hline

1 & 3^2{:}\SL_2(3) & 5^2{:}\SL_2(3) & \SL_2(3) & 3^2 & 5^2 & \ZZ_3 & 1\\

2 & 3^2{:}\GL_2(3) & 5^2{:}\GL_2(3) & \GL_2(3) & 3^2 & 5^2 & \S_3& \ZZ_2\\

3 & 3^2{:}\GL_2(3) & 7^2{:}\GL_2(3) & \GL_2(3) & 3^2 & 7^2 & \S_3 & 1\\

4 & 5^2{:}\GL_2(3) & 7^2{:}\GL_2(3) & \GL_2(3) & 5^2 & 7^2 & \ZZ_2 & 1 \\

\hline

5&\A_6&\PSL_2(11)&\A_5&6&11&\A_4&\D_6  \\

6&\PSL_2(11) &2^4{:}\A_5&\A_5&11&16&\D_6&\ZZ_2^2  \\

7&\A_6&2^4.\A_5&\A_5&6&16&\A_4&\ZZ_2^2  \\

8&\S_6&2^4.\S_5&\S_5&6&16&\S_4&\D_8  \\

9&\A_7&2^4.\A_6&\A_6&7&16&\A_5&\S_4  \\

10&\S_7&2^4.\S_6&\S_6&7&16&\S_5&2\times\S_4  \\

11&\A_8&2^4.\A_7&\A_7&8&16&\A_6&\PSL_2(7) \\ 

12&\A_9&2^4.\A_8&\A_8&9&16&\A_7&2^3{:}\GL_3(2) \\

13& \S_9  & \Sp_6(2) &\S_8 &  9 &36 & \S_7 & (\S_4\times\S_4).2\\ \hline

14&\A_7 & 2^3{:}\GL_3(2) & \PSL_2(7) &  15 & 8 & \A_4 & \S_4\\

15&5^2{:}\SL_2(5)&{11}^2{:}\SL_2(5)&\SL_2(5)&5^2&11^2&\ZZ_5&1\\

16&{13}^2{:}\SL_2(13)&3^6{:}\SL_2(13)&\SL_2(13)&13^2&3^6&\ZZ_{13}&\ZZ_3\\

\hline
\end{array}\]
\centerline{\bf Table~\ref{Intro}}

\vskip0.1in

\begin{corollary}\label{pair-2-trans}
Let $P$ and $Q$ be two $2$-transitive permutation groups which are not isomorphic but their stabilisers are isomorphic to the same group $H$.
Then $(P,Q,H)$ is one of the triples $(G_v,G_w,G_{uw})$ in Table~$\ref{Intro}$.
\end{corollary}

For locally 3-arc-transitive graphs, there are only six of such amalgams.

\begin{corollary}\label{faith-stab-3-trans}
Let $\Ga$ be a connected locally $(G,3)$-arc-transitive, and let $\{v,w\}$ be an edge of $\Ga$.
Assume that $G_v^{[1]}=G_w^{[1]}=1$.
Then $(G_v,G_w,G_{vw})$ is one of the following triples:
$(\A_7,\A_7,\A_6),\ (\S_7,\S_7,\S_6),\ (\A_7,2^4{:}\A_6,\A_6),\ (\S_7,2^4{:}\S_6,\S_6), (\A_8,2^4{:}\A_7,\A_7),$ $ (\A_9,2^4{:}\A_8,\A_8),$ or $\ (\S_9,\Sp_6(2),\S_8).$
\end{corollary}

\section{Examples}

Let $H$ be a 2-transitive permutation group of degree $n$, and let $p$ be a prime which is coprime to the order $|H|$.
Let $V=\FF_p^n$ be the permutation module of $H$ over $\FF_p$.
Then $H\leqslant\GL_n(p)$.
Let $\s$ be the unique involution in $\Z(\GL_n(p))$.
Then $\s$ reverses every vector in $V$, and $\l \s,H\r=\l \s\r\times H$ since $\s$ centralizes $H$ and $H$ is a 2-transitive permutation groups.
Define a group .
\[X=\ZZ_p^n{:}(\l \s\r\times H).\]
Let $a_1,a_2,\dots,a_n$ be a basis for $V=\ZZ_p^n$. Then $a_i^\s=a_i$, and $H$ is 2-transitive on the basis.
Let 
\[g=a_1\s.\]
Then $g$ is an involution.
Let
\[S=\{a_1\s,a_2\s,\dots,a_n\s\},\]
and let $R=\l S\r$.
Then $(a_i\s)(a_j\s)=a_ia_j\in R$, and $H$ acts 2-transitively on $S$.
Now $H$ induces a subgroup of $\Aut(R)$.
Let $\Ga=\Cay(R,S)$, and 
\[G=R{:}H.\]
Thus $\Ga$ is a $(G,2)$-arc-transitive Cayley graph of $R$ with vertex stabiliser $H$.


\begin{proposition}
Any $2$-transitive permutation group is the vertex stabiliser of $2$-arc-transitive graph.
\end{proposition}

\begin{example}
{\rm
Let $(G_v,G_w,G_{vw})$ be one of the triples listed in Table~\ref{Intro} except for 
rows~5, 6, 10, 12 and 13.
Then there exists a group $G$ such that $G=G_vG_w$, and thus there is a complete bipartite graph
$\Ga=\K_{m,n}$ which is locally $(G,2)$-arc-transitive.
These examples given in the next table.
}
\end{example}

%
%
%
%
%

\[\begin{array}{|c ccc c |}\hline
G & G_v&G_w&G_{vw}& \Ga \\ \hline

(3^2\times 5^2){:}\SL_2(3)& 3^2{:}\SL_2(3) & 5^2{:}\SL_2(3) & \SL_2(3) & \K_{3^2,5^2}\\

(3^2\times 5^2){:}\GL_2(3)& 3^2{:}\GL_2(3) & 5^2{:}\GL_2(3) & \GL_2(3) & \K_{3^2,5^2}\\

(3^2\times 7^2){:}\GL_2(3) & 3^2{:}\GL_2(3) & 7^2{:}\GL_2(3) & \GL_2(3) & \K_{3^2,7^2}\\

(5^2\times 7^2){:}\GL_2(3) & 5^2{:}\GL_2(3) & 7^2{:}\GL_2(3) & \GL_2(3) & \K_{5^2,7^2} \\

2^4{:}\A_6&\A_6&2^4.\A_5&\A_5&\K_{6,2^4}  \\

2^4{:}\S_6&\S_6&2^4.\S_5&\S_5&\K_{6,2^4}  \\

2^4{:}\A_7&\A_7&2^4.\A_6&\A_6&\K_{7,2^4}  \\

2^4{:}\A_8&\A_8&2^4.\A_7&\A_7&\K_{8,2^4}  \\

\A_8 &\A_7 & 2^3{:}\GL_3(2) & \GL_3(2) &  \K_{7,15}\\

(5^2\times11^2){:}\SL_2(5)&5^2{:}\SL_2(5)&{11}^2{:}\SL_2(5)&\SL_2(5)&\K_{5^2,11^2}\\

(13^2\times 3^6){:}\SL_2(13)&{13}^2{:}\SL_2(13)&3^6{:}\SL_2(13)&\SL_2(13)&\K_{13^2,3^6}\\

\hline
\end{array}\]

\vskip0.1in
\noindent {\bf Problem.} Construct explicit examples of graphs with the amalgams in row~5, 6, 10, 12, or 13
in Table~\ref{Intro}.
\vskip0.1in

%
%

\begin{example}
{\rm
It is well-known that the triple $(\A_7,\A_7,\A_6)$ is the amalgam for both the complete graph $\K_8$ and 
Hoffman-Singleton graph.
The former is 2-arc-transitive but not 3-arc-transitive, but the latter is 3-arc-transitive.
The same phenomena appears for $(\S_7,\S_7,\S_6)$.
}
\end{example}

\section{Composition factors of the stabilisers}\label{B-pty}

Let $\Ga=(V,E)$ be a connected $G$-edge-transitive bipartite graph with biparts $U$ and $W$.
For vertices $v_0,v_1,\dots,v_\ell$ of $\Ga$, let
\[G_{v_0v_1\dots v_\ell}^{[1]}=G_{v_0}^{[1]}\cap G_{v_1}^{[1]}\cap\dots\cap G_{v_\ell}^{[1]}.\]
Then, for an edge $\{v,w\}$, the group $G_{vw}^{[1]}$ is the point-wise stabiliser of the subset of vertices $\Ga(v)\cup\Ga(w)$.

\begin{lemma}\label{chain-1}
Let $v_0,v_1,\dots,v_\ell$ be an $\ell$-arc, where $\ell$ is a positive integer.
Then we have a chain of normal subgroups:
\[G_{v_0v_1\dots v_\ell}^{[1]}\lhd G_{v_0v_1\dots v_{\ell-1}}^{[1]}
  \lhd\dots\lhd G_{v_0v_1}^{[1]}\lhd G_{v_0}^{[1]}\lhd G_{v_0}.\]
\end{lemma}
\proof
For each integer $i$ with $0\leqslant i\leqslant\ell$, the group $G_{v_0v_1\dots v_i}^{[1]}$ is the kernel of
$G_{v_0v_1\dots v_i}$ acting on $\Ga(v_0)\cup\Ga(v_1)\cup\dots\cup\Ga(v_i)$, and hence
\[G_{v_0v_1\dots v_i}^{[1]}\lhd G_{v_0v_1\dots v_i}.\]
For $i\geqslant 1$, we have the inclusion relation 
$G_{v_0v_1\dots v_i}^{[1]}\leqslant G_{v_0v_1\dots v_{i-1}}^{[1]}\leqslant G_{v_0v_1\dots v_i}$,
and we thus conclude that $G_{v_0v_1\dots v_i}^{[1]}\lhd G_{v_0v_1\dots v_{i-1}}^{[1]}$.
This holds for all possible values of $i$ with $1\leqslant i\leqslant\ell$, giving the chain of normal subgroups as claimed.
\qed

\begin{lemma}\label{chain-2}
Let $\{v,w\}$ be an edge, and $v_0,v_1,\dots,v_\ell$ be an $\ell$-arc, where $\ell$ is a positive integer.
Then for $0\leqslant i\leqslant \ell-1$, the factor group 
$G_{v_0v_1\dots v_i}^{[1]}/G_{v_0v_1\dots v_{i+1}}^{[1]}$ is isomorphic to a subnormal subgroup of $G_{vw}^{\Ga(v)}$ or
$G_{vw}^{\Ga(w)}$.
\end{lemma}
\proof
Observing that $G_{v_0v_1\dots v_i}^{[1]}\lhd\lhd G_{v_i}^{[1]}\lhd G_{v_iv_{i+1}}$, we obtain
$$G_{v_0v_1\dots v_i}^{[1]}/G_{v_0v_1\dots v_{i+1}}^{[1]}
\cong (G_{v_0v_1\dots v_i}^{[1]})^{\Ga(v_{i+1})}
\lhd\lhd G_{v_iv_{i+1}}^{\Ga(v_{i+1})}.$$
Since $G$ is edge-transitive on $\Ga$, we have $G_{v_iv_{i+1}}\cong G_{vw}$,
and so
\[\mbox{ $G_{v_iv_{i+1}}^{\Ga(v_{i+1})}\cong G_{vw}^{\Ga(w)}$ or $G_{vw}^{\Ga(v)}$.}\]
Therefore, 
$G_{v_0v_1\dots v_i}^{[1]}/G_{v_0v_1\dots v_{i+1}}^{[1]}$ is isomorphic to a subnormal subgroup of $G_{vw}^{\Ga(v)}$ or
$G_{vw}^{\Ga(w)}$.
\qed

\vskip0.1in
{\bf Proof of Theorem~\ref{key-lem}:}
First of all, we claim that a  composition factor of $G_{v}^{[1]}$ is isomorphic to a composition factor of $G_{vw}^{\Ga(v)}$ or $G_{vw}^{\Ga(w)}.$ Let $T$ be a composition factor of  $G_v^{[1]}$.
Since $G\leqslant\Aut\Ga$ and $\Ga$ is finite and connected, there exists a path $v_0=v,v_1,\dots,v_\ell$ such that $G_{v_0v_1\dots v_\ell}^{[1]}=1$.
By Lemma~\ref{chain-1}, we have a chain of normal subgroups:
\[1=G_{v_0v_1\dots v_\ell}^{[1]}\lhd G_{v_0v_1\dots v_{\ell-1}}^{[1]}
  \lhd\dots\lhd G_{v_0v_1}^{[1]}\lhd G_{v_0}^{[1]}=G_v^{[1]}.\]

Let $i$ be the largest integer with $0\leqslant i\leqslant\ell$ such that $T$ is a composition factor 
of $G_{v_0v_1\dots v_i}^{[1]}$.
Then $T$ is not a composition factor of $G_{v_0v_1\dots v_{i+1}}^{[1]}$.
Hence $T$ is a composition factor of  the quotient group
$G_{v_0v_1\dots v_i}^{[1]}/G_{v_0v_1\dots v_{i+1}}^{[1]}.$
By Lemma~\ref{chain-2}, $G_{v_0v_1\dots v_i}^{[1]}/G_{v_0v_1\dots v_{i+1}}^{[1]}$ is isomorphic to 
a subnormal subgroup of $G_{vw}^{\Ga(w)}$ or $G_{vw}^{\Ga(v)}$, and so 
$T$ is a composition factor of $G_{vw}^{\Ga(w)}$ or $G_{vw}^{\Ga(v)}$ as we claimed.

Now let $T_1$ be a composition factor of $G_v$.
If $T_1$ is a composition factor of $G_v/G_v^{[1]}\cong G_v^{\Ga(v)}$, then we are done.
We thus assume that $T_1$ is not a composition factor of $G_v^{\Ga(v)}\cong G_v/G_v^{[1]}$.
Then $T_1$ is a composition factor of $G_v^{[1]}$, and so $T_1$ is a composition factor of $G_{vw}^{\Ga(w)}$ or $G_{vw}^{\Ga(v)}$. This proves part~(i).

Similarly, let $T_2$ be a composition factor of $G_{vw}$. If $T_2$ is a composition factor of $G_{vw}/G_v^{[1]}\cong G_{vw}^{\Ga(v)}$, then we are done. Now we assume that  $T_2$ is a composition factor of $G_v^{[1]}$, and so $T_1$ is a composition factor of $G_{vw}^{\Ga(w)}$ or $G_{vw}^{\Ga(v)}$. Thus part~(ii) of Theorem~\ref{key-lem} hold.

%
%
%

\begin{corollary}\label{key-lem-2}
Let $\Ga$ be a connected graph, and let $G\le\Aut\Ga$ be transitive on the edge set.
Let $\{v,w\}$ be an edge of $\Ga$.
Then the following statements hold:
\begin{itemize}
\item[(i)]  The vertex stabiliser $G_v$ is soluble if and only if $G_v^{\Ga(v)}$ and $G_{vw}^{\Ga(w)}$ are both soluble.

\item[(ii)]  $G_{vw}$ is soluble if and only if $G_{vw}^{\Ga(v)}$ and $G_{vw}^{\Ga(w)}$ are both soluble; in particular.

\item[(iii)] $G_v$, $G_w$ are both soluble if and only if $G_v^{\Ga(v)}$ and $G_w^{\Ga(w)}$ are both soluble.

\item[(iv)] $G_v^{[1]}$ and $G_w^{[1]}$ are both soluble if and only if 
$(G_v^{[1]})^{\Ga(w)}$ and $(G_w^{[1]})^{\Ga(v)}$  are both soluble.

\end{itemize}
\end{corollary}
\proof
Part~(i), and part~(ii)  follow from Theorem~\ref{key-lem}\,(1), and (2), respectively.

For part~(iii), it is obvious that if $G_v$, $G_w$ are both soluble then$G_v^{\Ga(v)}$ and $G_w^{\Ga(w)}$ are both soluble. So it is sufficient to show if both $G_v^{\Ga(v)}$ and $G_w^{\Ga(w)}$ is soluble, then both $G_v$ and $G_w$ are soluble.
Suppose that $G_v$ is insoluble. 
Let $T$ be an insoluble composition factor of $G_v$.
By Theorem~\ref{key-lem}~(1), $T$ is a composition factor of $G_v^{\Ga(v)}$ or $G_{vw}^{\Ga(w)}$.
Similarly, if $G_w$ is insoluble, then either $G_w^{\Ga(w)}$ or $G_{vw}^{\Ga(v)}$ is insoluble.
It implies that if both $G_v^{\Ga(v)}$ and $G_w^{\Ga(w)}$ is soluble, then both $G_v$ and $G_w$ are soluble.
This proves part~(iii).

Similarly, to prove part~{iv}, we need only to prove if both $(G_v^{[1]})^{\Ga(w)}$ and $(G_w^{[1]})^{\Ga(v)}$  are soluble 
then so are $G_v^{[1]}$ and $G_w^{[1]}$.
Suppose that $G_v^{[1]}$ is insoluble.
Arguing as in the second paragraph with $G_v^{[1]}$ in the place of $G_v$ shows that 
either $(G_v^{[1]})^{\Ga(w)}$ or $(G_{vw}^{[1]})^{\Ga(v)}$ is insoluble.
Similarly, if $G_w^{[1]}$ is insoluble then either $(G_w^{[1]})^{\Ga(v)}$ or $(G_{vw}^{[1]})^{\Ga(w)}$ is insoluble.
It implies that if both $(G_v^{[1]})^{\Ga(w)}$ and $(G_w^{[1]})^{\Ga(v)}$  are soluble, 
then so are $G_v^{[1]}$ and $G_w^{[1]}$.
\qed

The next lemma refines Lemma~\ref{chain-2}.

\begin{lemma}\label{chain-3}
Let $\{u,w\}$ be an edge with $u\in U$ and $w\in W$, and let 
\[u_0,w_0,\dots,w_{i-1},u_i,w_i,u_{i+1},\dots\] 
be a path
such that $u_0=u$ and $w_0=w$.
Then for $i\geqslant 0$, the following  are true: 
\begin{itemize}
\item[(i)] $G_{u_0\dots u_i}^{[1]}/G_{u_0\dots u_iw_i}^{[1]}$ is isomorphic to a subnormal subgroup of $G_{uw}^{\Ga(w)}$.

\item[(ii)] $G_{u_0\dots w_i}^{[1]}/G_{u_0\dots w_iu_{i+1}}^{[1]}$ is isomorphic to a subnormal subgroup of $G_{uw}^{\Ga(u)}$.
\end{itemize}
\end{lemma}
\proof
We notice that 
\[\begin{array}{l}
G_{u_0\dots u_i}^{[1]}/G_{u_0\dots u_iw_i}^{[1]}\cong (G_{u_0\dots u_i}^{[1]})^{\Ga(w_i)}\lhd\lhd G_{u_iw_i}^{\Ga(w_i)}\cong G_{uw}^{\Ga(w)},\\

G_{u_0\dots w_i}^{[1]}/G_{u_0\dots w_iu_{i+1}}^{[1]}\cong (G_{u_0\dots w_i}^{[1]})^{\Ga(u_{i+1})}\lhd\lhd G_{w_iu_{i+1}}^{\Ga(u_{i+1})}\cong G_{uw}^{\Ga(u)}.
\end{array}\]
Then the proof follows.
\qed

\section{2-Transitive permutation groups}\label{2-trans}

The proof of Theorem~\ref{faith-stab} depends on the structural properties of stabilisers of 
2-transitive permutation groups.
We state the point stabilisers in the following two theorems, refer to \cite[Tables~7.3 and 7.4]{Cameron-book} ,
and \cite{LiSS}.

The {\it socle} of a group $X$ is the subgroup of $X$ generated by all
minimal normal subgroups of $X$, denoted by $\soc(X)$.

\begin{theorem}\label{AS-2-trans}
Let $X$ be an almost simple $2$-transitive permutation group
on $\Omega$ of degree $k$, and let $T=\soc(X)$ and $\o,\o'\in\Omega$.
Then $X/T\cong X_\o/T_\o\cong X_{\o\o}/T_{\o\o'}\leqslant \calO\leqslant \Out(T)$,
and further, they satisfy the following table:
\end{theorem}
{\small
\[\begin{array}{|c|c|c|c|c|}\hline
\mbox{row} & T & T_\o & T_{\o\o'} & k\\ \hline 

1& \A_n,n\geq 5 & \A_{n-1} & \A_{n-2} & n  \\ \hline

2& \PSL_n(q) & q^{n-1}{:}(\SL_{n-1}(q).{q-1\over(n,q-1)}) & [q^{2(n-2)}]{:}(\GL_{n-2}(q).{q-1\over(n,q-1)})& {q^n-1\over q-1} 
 \\ \hline

3& \Sz(q), q>2 & [q^2]{:}\ZZ_{q-1} & \ZZ_{q-1} & q^2+1 \\ 

4& \Ree(q),q>3 & [q^3]{:}\ZZ_{q-1} & \ZZ_{q-1} & q^3+1 \\

5& \PSU_3(q),q\geq 3 & [q^3]{:}\ZZ_{(q^2-1)/(3,q+1)} 
 & \ZZ_{(q^2-1)/(3,q+1)} & q^3+1 \\ \hline

6& \Sp_{2n}(2), & \POmega_{2n}^-(2).2 & 2^{2n-2}{:}\POmega_{2n-2}^-(2).2
	& 2^{2n-1}-2^{n-1} \\ 

7& n\geqslant3 & \POmega_{2n}^+(2).2 & 2^{2n-2}{:}\POmega_{2n-2}^+(2).2
	& 2^{2n-1}+2^{n-1} \\ \hline

8& \PSL_2(11) & \A_5 & \S_3 & 11 \\ 

9& \M_{11} & \M_{10} & 3^2{:}\Q_8 & 11 \\ 

10& \M_{11} & \PSL_2(11) & \A_5 & 12 \\ 

11& \M_{12} & \M_{11} & \M_{10} & 12 \\
 
12& \A_7 & \PSL_2(7) & \A_4 & 15 \\ 

13& \M_{22} & \PSL_3(4) & 2^4{:}\A_5 & 22 \\ 

14& \M_{23} & \M_{22} & \PSL_3(4) & 23 \\ 

15& \M_{24} & \M_{23} & \M_{22} & 24 \\ 

16& \PSL_2(8) & \D_{18} & \ZZ_2 & 28 \\ 

17& \HS & \PSU_3(5).2 & \A_6.2^2 & 176 \\ 

18& \Co_3 & \McL{:}2 & \PSU_4(3){:}2 & 276 \\ \hline

\end{array}\]

\nobreak
{\small \centerline{{\bf Table~\ref{2-trans}.1}: 
Almost simple 2-transitive affine permutation groups and their stabilisers}}
}

\vskip0.08in
A complete list of finite 2-transitive permutation groups was obtained
by Hering, see for example \cite{Liebeck}.
Here we list their point stabilisers.

\begin{theorem}\label{HA-2-trans}
Let $X=\ZZ_p^e{:}X_\o$ be an affine $2$-transitive permutation group of degree
$p^e$ on a set $\Omega$, where $\o\in\Omega$.
Then for a point $\o'\in\Omega\setminus\{\o\}$, either 
\begin{enumerate} 
\item[(1)] $X_\omega\leqslant \GammaL_1(p^e)\cong\ZZ_{p^e-1}{:}\ZZ_e$, and 
$X_{\o\o'}\lesssim\ZZ_e$, or 

\item[(2)] $X_\o$ and $X_{\o\o'}$ are as in Table~$\ref{2-trans}.2$,
\end{enumerate}
\end{theorem}
{\small
\[\begin{array}{|c|c|c|c|c|c| c| }\hline
\mbox{\rm row} & X_\o & X_{\o\o'} & p^e & \\ \hline 

1 & \SL_n(q).\calO & q^{n-1}{:}\GL_{n-1}(q).\calO & q^n &  \calO\leqslant\ZZ_{q-1}{:}\ZZ_e\\

2 & \Sp_{2n}(q).\calO & q.q^{2n-2}{:}(\GL_1(q)\times \Sp_{2n-2}(q)).\calO &q^{2n} &  \calO\leqslant \ZZ_e\\

3 & \G_2(q).\calO & [q^5]{:}\GL_2(q).\calO & q^6 & \calO\leqslant \ZZ_{q-1}{:}\ZZ_{(3,q)e}\\  \hline

4& \Q_8{:}3, \ \Q_8.6,\ (\Q_8{:}3).4 & 1,\ \ZZ_2, [4]& 5^2& \\

5 &\Q_8{:}\S_3,\ 3\times(\Q_8.2),\ 3\times(\Q_8{:}\S_3) &  1,\ 1, \ \ZZ_3 & 7^2&\\

6&5\times(\Q_8{:}3), \ 5\times(\Q_8{:}\S_3) & 1, \ \ZZ_2 & 11^2&\\

7&{11}\times(\Q_8{:}\S_3) & 1 & 23^2&\\ \hline

8 &2^{1+4}{:}5,2^{1+4}{:}\D_{10},2^{1+4}{:}(5{:}4),& 2, [4], [8],& 3^4 &\\
 &2^{1+4}{:}\A_5,2^{1+4}{:}\S_5 &2.\A_4,2.\S_4&& \\ 
 \hline

9&\SL_2(13) & \ZZ_3 & 3^6 &\\ 

10&\A_6 & \S_4 & 2^4  &\\ 

11&\A_7 & \PSL_2(7) & 2^4 & \\ 

12&\PSU_3(3) & 4.\S_4,\ 4^2.\S_3  & 2^6 &\\ \hline

13&\SL_2(5),\ 5\times\SL_2(5) & 1,\ \ZZ_5 & 11^2 &\\

14& 9\times\SL_2(5) & \ZZ_3 &19^2&\\ 

15&7\times\SL_2(5),\ 7\times(4\circ\SL_2(5))  & 1,\ \ZZ_2 &29^2 &\\

16&{29}\times\SL_2(5) & 1 & 59^2 & \\ \hline

\end{array}\]
}
{\small
\centerline{{\bf Table~\ref{2-trans}.2}: Affine 2-transitive groups and their stabilisers} 
  }

\begin{lemma}\label{2-pt-stab}
 Let $X$ be a $2$-transitive group on $\Ome$ of degree $k$, and take $\o,\o'\in\Ome$.
 Then the stabiliser $X_{\o\o'}$ has a transitive representation of degree $k-1$ if and only if $X=\A_7$ or $\S_7$,
 and $X_{\o\o'}=\A_5$ or $\S_5$, respectively.
 \end{lemma}
 \proof
 Suppose $X_{\o\o'}$ has has a transitive permutation representation of degree $k-1$. Then $k-1$ divides $|X_{\o\o'}|$. It is clear that $k-1$ does not divide $|X_{\o\o'}|$ for candidates in rows 3 to 5 and 8 to 18 in Table~\ref{AS-2-trans} and rows 4 to 16 in Table~\ref{HA-2-trans}. 
 
 Now we inspect the other candidates given  in Tables~\ref{AS-2-trans} and \ref{HA-2-trans}.
 Suppose $\soc(X)=\A_k$ as in row 1 of Table~\ref{AS-2-trans}. Then $X_{\o\o'}=\A_{k-2}$ or $\S_{k-2}$.
If $X_{\o\o'}=\A_{k-2}$ or $\S_{k-2}$ has a transitive permutation representation of degree $k-1$, then $k=7$.

Suppose that $\soc(X)=\PSL_n(q)$, with $k={q^n-1\over q-1}$ as in row 2 of Table~\ref{AS-2-trans}. Suppose further $X_{\o\o'}$ has a transitive permutation representation of degree $k-1=\frac{q(q^{n-1}-1)}{q-1}$. If $q^{n-1}-1$has a primitive prime divisor $r$ say, then $r$ does not divides the order $|X_{\o\o'}|=|[q^{2(n-2)}]{:}(\GL_{n-2}(q).\ZZ_{q-1\over(n,q-1)}).\calO|$, which is not possible. Thus $q^{n-1}-1$ has no primitive prime divisor, and $n-1=2$ or $(q,n-1)=(2,6)$ by Zsigmondy Theorem. For the former, $k-1=q(q+1)$, and $|X_{\o\o'}|=q^2(q-1)^2/(3,q-1)$ which is not divisible by $q+1$ So  $k-1=126$, and $X_{\o\o'}=[2^{10}]:GL(5,2)$. But  $X_{\o\o'}$ has no transitive permutation representation of degree 126 contradicts to our assumption.

Similarly, let $X=q^n{:}\SL_n(q).\calO$ be as in the first row of Table~\ref{HA-2-trans}.
If $k-1=q^n-1$ divides the order of $X_{\o\o'}=q^{n-1}{:}\GL_{n-1}(q).\calO$, then $q^n-1$ has no primitive prime divisor.
Thus, by Zsigmondy Theorem, $n=2$ or $(q,n)=(2,6)$.
For the former, $k-1=q^2-1$, and $|X_{\o\o'}|=q(q-1)$ which is not divided by $k-1$. For the latter, $k-1=63$, $X_{\o\o'}=[2^5]:\GL(4,2)$ which has no transitive permutation representation of degree 63.

For candidates in rows 6,7 in Table~~\ref{AS-2-trans}, and rows 2,3 in Table~\ref{HA-2-trans}, we have the following table:

{\small
\[\begin{array}{|c|c|c| }\hline

 X_\o & X_{\o\o'} & k-1  \\ \hline 
 
 \POmega_{2n}^-(2).2.\calO  & 2^{n^2-n+1}(2^{n-1}+1)\Pi_{i=1}^{n-2}(2^{2i}-1).\calO 
	& 2^{2n-1}-2^{n-1} \\ 

 \POmega_{2n}^+(2).2.\calO  & 2^{n^2-n+1}(2^{n-1}-1)\Pi_{i=1}^{n-2}(2^{2i}-1).\calO 
	& 2^{2n-1}+2^{n-1} \\ \hline


 \Sp_{2n}(q).\calO & q^{n^2}(q-1)\Pi_{i=1}^{n-1}(q^{2i}-1).\calO &q^{2n}-1 \\

 \G_2(q).\calO & q^6(q-1)^2(q+1).\calO & q^6 -1\\  \hline

\end{array}\]
}

It is easy to see that $k-1$ does not divide the order $|X_{\o\o'}|$. So $X_{\o\o'}$ has no transitive permutation representation of degree $k-1$.\qed

\section{Proof of Theorem~\ref{faith-stab}}\label{main-thm}

Let $\Ga$ be a connected locally $(G,2)$-arc-transitive graph, and let $\{v,w\}$ be an edge of $\Ga$. Assume that $G_v^{[1]}=G_w^{[1]}=1$. In this section, we will give the proof of Theorem~\ref{faith-stab}.
We first notice that, since $G_v\cong G_v^{\Ga(v)}$ and $G_w\cong G_w^{\Ga(w)}$, 
the stabilisers $G_v,G_w$ are 2-transitive permutation groups 
which share common  point stabilisers $G_{vw}$.

To prove this theorem, we will split it into several cases, given in the following few lemmas.

\begin{lemma}\label{isomorphic socle}
Suppose $G_v$ and $G_w$ have isomorphic socles. Then $G_v\cong G_w$ and $\Ga$ is regular.
\end{lemma}

\proof Since $G_v$ and $G_w$ are 2-arc-transitive groups, we have either they are both affine groups or they are both almost simple groups. Suppose $\soc(G_v)\cong\soc(G_w)\cong\ZZ_p^d$, then $G_v\cong \ZZ_p^d{:}G_{vw}\cong G_w$, and 
$\Ga$ is regular of valency $p^d$. Suppose that $G_v$ and $G_w$ are almost simple and have isomorphic socles.
Since they share point stabiliser, we have $G_v\cong G_w$ by checking the candidates in Table~\ref{2-trans}.1.\qed

From now on, we assume that $G_v$ and $G_w$ have non-isomorphic socles. We need to search for different 2-transitive groups which have the same point stabiliser, by inspecting 
the 2-transitive permutation groups, listed in Tables~\ref{2-trans}.1 and \ref{2-trans}.2.

\begin{lemma}\label{both soluble}
Suppose both $G_v$ and $G_w$ are soluble. Then 
\[\begin{array}{ l | l l l l l |} 
G_v & 3^2{:}\SL_2(3) & 3^2{:}\GL_2(3) & 3^2{:}\GL_2(3)& 5^2{:}\GL_2(3) \\ \hline

G_w & 5^2{:}\SL_2(3)& 5^2{:}\GL_2(3) & 7^2{:}\GL_2(3) & 7^2{:}\GL_2(3)\\
\end{array}\]
\centerline{Table \ref{main-thm}.1: Soluble stabilisers}

\end{lemma}

\proof Suppose that $G_v$ and $G_w$ are soluble.
Then by inspecting 
the 2-transitive permutation groups, listed in Tables~\ref{2-trans}.1 and \ref{2-trans}.2, we have 

\[\begin{array}{ ll } \hline
G_{vw} & |\Ga(v)|  \\ \hline

\leqslant \ZZ_{p^d-1}{:}\ZZ_d & p^d \\

\GL_2(2) & 2^2 \\

\SL_2(3),\ \GL_2(3) & 3^2 \\

\Q_8{:}3, \ \Q_8.6,\ (\Q_8{:}3).4  & 5^2 \\

\Q_8{:}\S_3,\ 3\times(\Q_8.2),\ 3\times(\Q_8{:}\S_3)  & 7^2 \\

5\times(\Q_8{:}3), \ 5\times(\Q_8{:}\S_3)  & 11^2 \\

{11}\times(\Q_8{:}\S_3) &  23^2 \\ 

2^{1+4}{:}5,2^{1+4}{:}\D_{10},2^{1+4}{:}(5{:}4) & 3^4 \\  \hline

\end{array}\]
\nobreak
\centerline{Table \ref{main-thm}.2}

\vskip0.08in
We note that the candidate in the second line, i.e.,  $G_{vw}=\GL_2(2)$, is the same as
the one in the first line with $p^d=2^2$ and $G_{vw}=\ZZ_{p^d-1}{:}\ZZ_2\cong\S_3$.

Assume first that $G_{vw}\leqslant\ZZ_{p^d-1}{:}\ZZ_d$ with $|\Ga(v)|=p^d$.
Then $G_{vw}$ is a split metacyclic group. So $G_{vw}$ can not appear as one  of rows~3-8 of Table \ref{main-thm}.2.

We thus have that $G_{vw}$ only appears in one of rows~3-8 of Table \ref{main-thm}.2.
It follows that  $G_v,G_w$ are as in Table \ref{main-thm}.1.\qed

\begin{lemma}\label{one insoluble}
Suppose $G_v$ is insoluble. Then $G_{vw}$ is insoluble.
\end{lemma}
\proof Suppose $G_{vw}$ is soluble. Then one of the candidates in the following table appears:
\[\begin{array}{ lll ll } \hline

G_v & G_{vw} & |\Ga(v)| & \calO & \\ \hline

\A_5, \S_5 & \A_4,\S_4 & 5 & & \\
\PSL_3(3).\calO & 3^2{:}\SL_2(3).\calO&26 &&\\
\PSL_2(q).\calO & [q]{:}\ZZ_{(q-1)/(2,q-1)}.\calO & q+1 & \calO\leqslant (2,q-1).f & q=p^f\geqslant4\\

\Sz(q).\calO & [q^2]{:}\ZZ_{q-1}.\calO & q^2+1 & \calO\leqslant f & q=2^f>2, f \mbox{ odd}\\

\Ree(q).\calO & [q^3]{:}\ZZ_{q-1}.\calO & q^3+1 & \calO\leqslant f & q=3^f>3, f \mbox{ odd} \\

\PSU_3(q).\calO & [q^3]{:}\ZZ_{(q^2-1)/(3,q+1)} & q^3+1  & \calO\leqslant (3,q+1).f & q\geq 3\\ 

\PGammaL_2(8) & 9{:}6 & 28 & &\\ \hline
\end{array}\]
\centerline{Table \ref{main-thm}.3: Soluble edges stabilisers}
\vskip0.08in

We note that, in Table~\ref{main-thm}.2,  the candidate in the first line is the same as the one in the second line with $q=4$ due to $\A_5\cong\PSL_2(4)$.
The rest  are mutually non-isomorphic.
Thus $G_w$ should be soluble, appearing in Table~\ref{main-thm}.1,
and we conclude that in this case $G_v$ and $G_w$ do not have isomorphic stabilisers. So $G_{vw}$ is insoluble.\qed

\begin{lemma}\label{$G_{vw}$ almost simple}
Suppose $G_{vw}$ is almost simple. Then $(G_v,G_w,G_{vw})$ are listed in rows 5 to 14 of Table~\ref{Intro}.
\end{lemma}

\proof
First, we observe that $G_{vw}$ does not have socle isomorphic to one of the following groups:
\begin{quote}
$\POmega_{2n}^-(2)$ with $n\geqslant3$, $\POmega_{2n}^+(2)$ with $n\geqslant4$,  $\M_{11}$, 
$\M_{22}$, $\M_{23}$, $\PSU_3(5)$, $\McL$, $\Sp_{2n}(q)$ with $(n,q)\not=(2,2)$, $\G_2(q)$, $\PSU_3(3)$,
\end{quote}
since each of them  appears only one time as the socle of the stabiliser of 2-transitive permutation groups. Similarly, $G_{vw}$ is not $\M_{10}$.
This shows that $\soc(G_{vw})$ is an alternating group of a linear group.

{\bf Case 1.} Let first $\soc(G_{vw})=\A_m$, where $m\geqslant5$.
Notice that the isomorphisms between  $\A_m$ and other simple groups:
\[\begin{array}{l}
\A_5\cong\PSL_2(5)\cong\PSL_2(4) \cong\POmega_4^-(2), \\ 

\A_6\cong\PSL_2(9)\cong\Sp_4(2)', \\

\A_8\cong\PSL_4(2)\cong\POmega_6^+(2).
\end{array}\]
Now $G_w$ is one of such groups such that $\soc(G_w)\not\cong\soc(G_v)$.
It implies that $m\leqslant8$.

Notice that $\PSL_2(4)\cong\SL_2(4)$, $\PSL_4(2)\cong\SL_4(2)$, and $\Sp_4(2)\cong \S_6$, we have  $G_v$ and $G_w$ is one of the following groups:

\[\begin{array}{c |c|c|c|c|c|c} 

G_v\mbox{:almost simple} & \A_6,\S_6 &\A_7,\S_7& \A_8,\S_8 &\S_9,\S_9& \Sp_4(2)&\PSL_2(11)\\ \hline

G_{vw} & \A_5,\S_5&\A_6,\S_6&\A_7,\S_7&\A_8,\S_8&\POmega_6^+(2)&\A_5\\ 
\end{array}\]
\vskip0.1in
\[\begin{array}{c |c|c|c|c|c} 

G_v\mbox{: affine}&4^2{:}\SL_2(4),4^2{:}\SL_2(4).2 &2^4{:}\SL_4(2)&2^4{:}\Sp_4(2),&2^4{:}\A_6&2^4{:}\A_7\\ \hline

G_{vw} &\SL_2(4),\SL_2(4).2&\SL_4(2)&\Sp_4(2)&\A_6&\A_7\\ 
\end{array}\]

\centerline{Table \ref{main-thm}.4: altnating edge stabilisers}

Suppose $m=5$. Then $G_v,G_w\rhd\A_6$, $\PSL_2(11)$, $2^4{:}\SL_2(4)$.
It then follows that $\{G_v,G_w\}$ is one of the following four pairs, noticing that
$\SL_2(4)\cong\A_5$ and $\SL_2(4).2\cong\SigmaL_2(4)\cong\S_5$:
\[\mbox{$\{\A_6,\PSL_2(11)\}$, $\{\PSL_2(11),2^4{:}\A_5\}$,
$\{\A_6,2^4{:}\A_5\}$, or $\{\S_6,2^4{:}\S_5\}$,}\]
as listed in rows~5-8 of Table~\ref{Intro}.

If $m=6$, then each of $ G_v$ and $G_w$ is conjugate to $\A_7$, $\S_7$ $2^4{:}\Sp_4(2)$, or $2^4{:}\A_6$, and hence
$\{G_v,G_w\}=\{\A_7,2^4{:}\A_6\}$ or $\{\S_7,2^4{:}\S_6\}$, 
as in rows~9-10 of Table~\ref{Intro}.

If $m=7$, then $G_v\rhd\A_8$ or $2^4{:}\A_7$, and so is $G_w$.
This  leads to 
$\{G_v,G_w\}=\{\A_8,2^4{:}\A_7\}$, 
as in row~11 of Table~\ref{Intro}.

If $m=8$, then each of $G_v$ and $G_w$ is $\A_9$, $\Sp_6(2)$, or $2^4{:}\GL_4(2)$, leading to 
$\{G_v,G_w\}=\{\A_9,2^4{:}\A_8\}$ or $\{\S_9,\Sp_6(2)\}$, 
giving in rows~12-13 of Table~\ref{Intro}.

\vskip0.1in
{\bf Case 2.} Assume that $\soc(G_{vw})$ is a linear group which is not isomorphic to $\A_m$ with $m\geqslant5$.
We notice that $\PSL_2(7)\cong\SL_3(2)\cong \GL_3(2)$.
Then $\soc(G_{vw})=\PSL_2(7)$, and $G_v,G_w\rhd\A_7$ or $2^3{:}\GL_3(2)$.
Therefore, we conclude that $\{G_v,G_w\}=\{\A_7,2^3{:}\GL_3(2)\}$, 
which is listed in row~14 of Table~\ref{Intro}.\qed

\begin{lemma}\label{$G_{vw}$ non-almost simple}
Assume $G_{vw}$ is insoluble but not almost simple. Then $(G_v,G_w,G_{vw})$ are listed in rows 15 and 16 of Table~\ref{Intro}.
\end{lemma}
\proof
Assume $G_{vw}$ is insoluble but not almost simple.
All the possibilities of $G_{vw}$ are listed in the following table.
\[\begin{array}{cc c}\hline
G_{vw} & \soc(G_v) & \\ \hline 

q^n{:}\SL_n(q).\calO & \PSL_{n+1}(q) &\\

\SL_n(q).\calO & q^n &\\

\Sp_{2n}(q).\calO & q^{2n} & \\  

2^{1+4}{:}\A_5,2^{1+4}{:}\S_5 & 2^4 \\  

\SL_2(13) & 3^6 \\ 

\SL_2(5),\ 5\times\SL_2(5) & 11^2 \\

9\times\SL_2(5) &19^2 \\ 

7\times\SL_2(5),\ 7\times(4\circ\SL_2(5))  &29^2  \\

{29}\times\SL_2(5) & 59^2 \\ \hline

\end{array}\]
The candidate in the first line or the third line is not isomorphic to any other one in the table.
This leads to the following possibilities:
\[\{G_v,G_w\}=\{5^2{:}\SL_2(5),11^2{:}\SL_2(5)\},\ \mbox{or}\ \{13^2{:}\SL_2(13), 3^6{:}\SL_2(13)\}.\]
These are listed in rows~15-16 of Table~\ref{Intro}.\qed

 {\bf Proof of Theorem~\ref{Intro}:}  By assumption, $G_{v}^{[1]}=G_w^{[1]}=1$, so  $G_v\cong G_v^{\Ga(v)}$, $G_w\cong G_w^{\Ga(w)}$, 
and the stabilisers $G_v,G_w$ are 2-transitive permutation groups 
which share common  point stabilisers $G_{vw}$. Suppose first that $G_v$ and $G_w$ have isomorphic soles. Then by Lemma~\ref{isomorphic socle}, $G_v\cong G_w$ and $\Ga$ is regular. 

Now suppose $G_v$ and $G_w$ have non-isomorphic socles. Suppose further that both of them are soluble. Then by Lemma~\ref{both soluble}, $(G_v,G_w,G_{vw})$ are listed in the first four rows of Table~\ref{Intro}. 
Suppose $G_v$ is insoluble. Then by Lemma~\ref{one insoluble}, $G_{vw}$ is insoluble. Thus by Lemmas~\ref{$G_{vw}$ almost simple} and \ref{$G_{vw}$ non-almost simple}, $(G_v,G_w,G_{vw})$ are listed in rows 5 to 16 of Table~\ref{Intro}. Thus the theorem holds.\qed

{\bf Proof of Corollary~\ref{faith-stab-3-trans}:}
We first notice that $\Ga$ is locally $(G,3)$-arc-transitive if and only if 
$G_{w_1vw}$ and $G_{vwv_1}$ are transitive on $[G_{w_{-1}vw}:G_{w_{-1}vwv_1}]$ and
$[G_{vwv_1}:G_{w_{-1}vwv_1}]$, respectively.
This is equivalent to $|G_{w_{-1}vw}:G_{w_{-1}vwv_1}|=|\Ga(w)|-1$ and
$|G_{vwv_1}:G_{w_{-1}vwv_1}|=|\Ga(v)|-1$.

Suppose first that $\Ga$ is a regular graph.
Then the 2-arc stabiliser $G_{vwv_1}$ is transitive on $\Ga(v)\setminus\{w\}$.
Observe that $G_{vwv_1}=(G_w)_{vv_1}$ is the stabiliser of the two points $v,v_1\in\Ga(w)$
in the 2-transitive permutation group $G_w$.
By Lemma~\ref{2-pt-stab}, we conclude that $(G_v,G_w,G_{vw})=(\A_7,\A_7,\A_6)$, or $(\S_7,\S_7,\S_6)$.

Now assume that $\Ga$ is not a regular graph.
We need analyse the triple of stabilisers $(G_v,G_w,G_{vw})$ listed in Table~\ref{Intro}.
We notice that, for each of the candidates in rows~1-8 and 14-16 of Table~\ref{Intro},
the order $|G_{vwv_1}|$ is not divisible by $|\Ga(v)|-1$.
Therefore, we conclude that $(G_v,G_w,G_{vw})$ can only be one of $(\A_7,2^4{:}\A_6,\A_6)$, $(\S_7,2^4{:}\S_6,\S_6)$, 
$(\A_8,2^4{:}\A_7,\A_7)$, $(\A_9,2^4{:}\A_8,\A_8)$, or $(\S_9,\Sp_6(2),\S_8)$.

For the first candidate $(G_v,G_w,G_{vw})=(\A_7,2^4{:}\A_6,\A_6)$, the intersection $G_{w_{-1}vwv_1}$ of
$G_{w_{-1}vw}=\A_5$ and $G_{vwv_1}=\S_4$ is a subgroup of $\ZZ_2^2$.
Since $|G_{w_{-1}vw}:G_{w_{-1}vwv_1}|\leqslant|\Ga(w)|-1=15$ and
$|G_{vwv_1}:G_{w_{-1}vwv_1}|=|\Ga(v)|-1=6$.
We thus have 
\[G_{w_{-1}vwv_1}=G_{w_{-1}vw}\cap G_{vwv_1}=\A_5\cap\S_4=\ZZ_2^2.\]
Therefore, $|G_{w_{-1}vw}:G_{w_{-1}vwv_1}|=15=|\Ga(w)|-1$ and
$|G_{vwv_1}:G_{w_{-1}vwv_1}|=6=|\Ga(v)|-1$, and $\Ga$ is locally $(G,3)$-arc-transitive.
Similarly, $(G_v,G_w,G_{vw})=(\S_7,2^4{:}\S_6,\S_6)$ is an amalgam for locally $(G,3)$-arc-transitive graph.

For the case $(G_v,G_w,G_{vw})=(\A_8,2^4{:}\A_7,\A_7)$, we have
\[G_{w_{-1}vw}=\A_6,\ G_{vwv_1}=\PSL_2(7).\]
Since  $|\Ga(v)|-1=8-1=7$, the stabiliser $G_{w_{-1}vwv_1}$ is of index at most 7 in $G_{vwv_1}=\PSL_2(7)$.
It implies that $|G_{vwv_1}:G_{w_{-1}vwv_1}|=7$, and $G_{w_{-1}vwv_1}=\S_4$.
Then the index $|G_{w_{-1}vw}:G_{w_{-1}vwv_1}|=|\A_6:\S_4|=15=|\Ga(w)|-1$.
We conclude that $G_{w_{-1}vw}$ is transitive on $\Ga(w)\setminus\{v\}$, and
$G_{vwv_1}$ is transitive on $\Ga(v)\setminus\{w\}$.
Therefore, $\Ga$ is locally $(G,3)$-arc-transitive.

For the triple $(G_v,G_w,G_{vw})=(\A_9,2^4{:}\A_8,\A_8)$,
the index of $G_{w_{-1}vwv_1}$ in $G_{w_{-1}vw}=\A_7$ is at most $|\Ga(w)|-1=15$ and
in $G_{vwv_1}=2^3{:}\GL_3(2)$ is at most $|\Ga(v)|-1=8$.
It follows that $G_{w_{-1}vwv_1}=\GL_3(2)$, and so  $|G_{w_{-1}vw}:G_{w_{-1}vwv_1}|=|\Ga(v)|-1$ and
$|G_{vwv_1}:G_{w_{-1}vwv_1}|=|\Ga(w)|-1$.
Therefore, $(\A_9,2^4{:}\A_8,\A_8)$ is a locally 3-arc-transitive amalgam.

For the case $(G_v,G_w,G_{vw})=(\S_9,\Sp_6(2),\S_8)$, we have
$G_{w_{-1}vw}=\S_7,\ G_{vwv_1}=(\S_4\times\S_4).2.$
The intersection $G_{w_{-1}vwv_1}$ of $G_{w_{-1}vw}$ and $G_{vwv_1}$ is $\S_4\times \S_3$ with index $8$ in $G_{vwv_1}$ and index $35$ in $G_{w_{-1}vw}$. Moreover, $|\Ga(v)|-1=9-1=8$ and $\Ga(w)|-1=35$.
We conclude that $G_{w_{-1}vw}$ is transitive on $\Ga(w)\setminus\{v\}$, and
$G_{vwv_1}$ is transitive on $\Ga(v)\setminus\{w\}$.
Therefore, $\Ga$ is locally $(G,3)$-arc-transitive.
\qed

\end{document}